\documentclass[12pt]{article}
\usepackage[reqno]{amsmath}
\usepackage{amssymb, theorem, enumerate}
\usepackage{graphicx}
\usepackage{tikz}
\usetikzlibrary{decorations.markings}
\usepackage{array}
\usepackage{float}
\usepackage{wrapfig}
\usepackage{fullpage,url}

\expandafter\ifx\csname pdfpagebox\endcsname\relax\else
\fi

\theoremstyle{change}
{\theorembodyfont{\slshape}
\newtheorem{theorem}{Theorem.}[section]
\newtheorem{lemma}[theorem]{Lemma.}
}

\theorembodyfont{\rmfamily}

\newcommand\cref[1]{Corollary~\ref{cor:#1}}

\def\proof{\noindent{{\sl Proof. }}}
\def\sqr#1#2{{\vbox{\hrule height.#2pt
    \hbox{\vrule width.#2pt height#1pt \kern#1pt
        \vrule width.#2pt}\hrule height.#2pt}}}
\def\eqed{\sqr53}
\def\qed{%
    \ifmmode\eqno\eqed
    \else\nobreak\ \hfill\eqed\medbreak\fi}

\newcommand\al{\alpha}
\newcommand\be{\beta}

\newcommand\ga{\gamma}

\newcommand\la{\lambda}

\newcommand\Om{\Omega}

\newcommand\cD{{\mathcal D}}

\newcommand\cF{{\mathcal F}}

\newcommand\cW{{\mathcal W}}

\newcommand\comp[1]{{\mkern2mu\overline{\mkern-2mu#1}}}
\newcommand\sbs{\subseteq}

\newcommand\seq[3]{#1_{#2},\ldots,#1_{#3}}

\newcommand\ip[2]{\langle#1,#2\rangle}

\DeclareMathOperator{\tr}{tr}

\newcommand\bW{\comp{W}}

\DeclareMathOperator\elsm{sum}
\newcommand\nw{N_{\cW}}

\title{Using the existence of t-designs to prove Erd\H{o}s-Ko-Rado}
\author{Chris Godsil\footnote{University of Waterloo, Waterloo, Canada. email: \protect\url{cgodsil@uwaterloo.ca}. C. Godsil gratefully acknowledges the support of the Natural Sciences and Engineering Council of Canada (NSERC), Grant No. RGPIN-9439.}, Krystal Guo\footnote{Universit\'{e} libre de Bruxelles, Brussels, Belgium. Part of this research was done when K. Guo was a post-doctoral fellow at University of Waterloo. email: \protect\url{krystal.guo@ulb.ac.be}}}

\begin{document}
\maketitle

\begin{abstract}
In 1984, Wilson proved the Erd\H{o}s-Ko-Rado theorem for $t$-intersecting families of $k$-subsets
of an $n$-set: he showed that if $n\ge(t+1)(k-t+1)$ and $\cF$ is a family of $k$-subsets
of an $n$-set such that any two members of $\cF$ have at least $t$ elements in common, then
$|\cF|\le\binom{n-t}{k-t}$. His proof made essential use of a matrix whose origin is not obvious.
In this paper we show that this matrix can be derived, in a sense, as a projection of
$t$-$(n,k,1)$ design. 
\end{abstract}

\section{Introduction}\label{sec:intro}

A family of sets is $t$-intersecting if every two sets in the family have at least $t$
elements in common. The Erd\H{o}s-Ko-Rado theorem states that if $\cF$ a $t$-intersecting family of
sets of size $k$ chosen from a set $N$ of size $n$ and $n\ge (t+1)(k-t+1)$, then 
\[
	|\cF| \ge\binom{n-t}{k-t}.
\]
If $n>(t+1)(k-t+1)$, equality holds if and only if $\cF$ consists of the $k$-subsets that contain a given
set of $t$ points from $V$. The lower bound on $n$ is necessary, because the result is false 
when the bound fails. Subsequently Ahlswede and Khachatrian \cite{AhlKha96, AhlKha97} determined the maximal families 
for all $n$. The result as just stated was proved by Wilson in 1984 \cite{Wilson1984}. 

The goal of this paper is to motivate a key step in Wilson's proof. He introduces a ``magic matrix''
with rows and columns indexed by the $k$-subsets of a $v$-set; he then determines the eigenvalues
of this matrix and, given these, fairly standard machinery then leads to the proof of the EKR-bound.
From private discussions with Rick Wilson, it is clear that this matrix was the result of
a lot of calculation and a lot of inspiration. Our aim in this paper is to present a derivation
which requires less effort and less brilliance. To this end, we give a simpler formulation of this matrix and
 show that it is equivalent to that of Wilson, using the recent proof of the existence of $t$-designs of Keevash \cite{ Kee14}.

\section{The Johnson Scheme}

Assume $N=\{1,\ldots,n\}$. The \textsl{Johnson scheme} $J(n,k)$ is a set of $01$-matrices
$\seq A0k$, with rows and columns indexed by the $k$-subsets of $N$, where $(A_r)_{\al,\be}=1$
if $|\al\cap\be|=k-r$ for $r = 0,\ldots,k$. We see that $A_0=I$. The matrices $\seq A1k$ are adjacency matrices
of graphs $\seq X1k$, where $X_1$ is the so-called \textsl{Johnson graph}. It can be shown that
two $k$-subsets are adjacent in $X_r$ if and only if they are at distance $k-r$ in the Johnson
graph. The Johnson scheme is discussed in detail in \cite[Chapter~6]{cg-km}, and anything we 
state here without proof is treated there.

The matrices $A_r$ satisfy
\[
	\sum_r A_r = J.
\]
Further, there are scalars $p_{i,j}(r)$ such that, for all $i$ and $j$,
\[
	A_iA_j = \sum_r p_{i,j}(r) A_r.
\]
Since the product of two symmetric matrices is symmetric if and only if the matrices commute, it
follows that
the space of the matrices $A_r$ is a commutative matrix algebra. (To use the standard jargon,
the matrices $\seq A0k$ form a \textsl{symmetric association scheme}, and their span is known
as the \textsl{Bose-Mesner algebra} of the scheme.) All matrices that occur in Wilson's proof
lie in the Bose-Mesner algebra of the Johnson scheme.

To define his matrix, Wilson used another basis for the Bose-Mesner algebra of the Johnson scheme.
Let $W_{i,j}(n)$ denote the matrix with rows indexed by the $i$-subsets of $N$, columns
indexed by the $j$-subsets of $N$ and with $(\al,\be)$-entry equal to 1 if $\al\sbs\be$.
(So each row of $W_{i,j}(n)$ sums to $\binom{n-i}{j-i}$.)
Let $\bW_{i,j}(n)$ denote the matrix with rows indexed by the $i$-subsets of $N$, columns
indexed by the $j$-subsets of $N$ and with $(\al,\be)$-entry equal to 1 if $\al\cap\be=\emptyset$.
Now define matrices $\seq D0k$ by
\[
	D_i = W_{i,k} \comp{W}_{i,k}^T
\]
(For details concerning these matrices see Wilson's paper \cite{}, or \cite[Section 6.4]{cg-km}.
Despite appeafances, these matrices are symmetric.)
The matrix $\Omega(n,k,t)$ is given by
\[
	\Omega(n,k,t) = \sum_{i=0}^{k-i} (-1)^{t-1-i} 
		\frac{\binom{k-1-i}{k -t}}{\binom{n-k-t+1}{k-t}} D_{k-i}.
\]
The matrices $I+\Om(n,k,t)$ form the key to Wilson's proof of the EKR theorem. We define
\[
	\nw(n,k,t) = I + \Om(n,k,t)
\]
and abbreviate $\nw(n,k,t)$ to $\nw$ where possible.

We use $M\circ N$ to denote the \textsl{Schur product} of two matrices of the same order, thus
\[
	(M\circ N)_{i,j} = M_{i,j}N_{i,j}.
\]
Since the set $\{0,\seq A0k\}$ is closed inder the Schur product, it follows that the Bose-Mesner
algebra of the Johnson scheme is closed under Schur product. 

The pertinent properties of $\nw$ are summarized in the following:

\begin{theorem}
	The matrix $\nw(n,k,t)$ is positive semidefinite and lies in the span of the matrices 
	$\seq A{k-t+1}{k}$.\qed
\end{theorem}

Wilson's proof that $\nw$ is positive semidefinite is highly non-trivial; it is presented
at somewhat greater length, but with no essential improvement, in \cite[Chapter~8]{cg-km}.

\section{Projections on to matrix algebras}

We use $\elsm(M)$ to denote the sum of the entries of a matrix $M$. We note that
\[
	\tr(M^TN) = \elsm(M\circ N)
\]
and so we have two expressions for the standard inner product on real matrices:
\[
	\ip MN = \tr(M^TN) = \elsm(M\circ N).
\]
Relatve to this inner product, the Schur idempotents $\seq A0k$ form an orthogonal basis for the Bose-Mesner algebra. We also observe that
\[
	\tr(M) = \ip{I}M,\quad \elsm(M) = \ip{J}M.
\]

We state a version of a result known as the clique-coclique bound. It is proved, for general
association schemes, as Lemma~3.8.1 in \cite{cg-km}. 

\begin{lemma}
	\label{lem:clqcoclq}
	Assume $v=\binom{n}{k}$. If $M$ and $N$ are matrices in the Bose-Mesner algebra of the Johnson
	scheme and
	\begin{enumerate}[(a)]
		\item 
		$M$ and $N$ are positive semidefinite, and
		\item
		for some constant $\ga$ we have $M\circ N=\ga I$,
	\end{enumerate}
	then
	\[
		\frac{\elsm(M)}{\tr(M)} \frac{\elsm(N)}{\tr(N)} \le v.\qed
	\]
\end{lemma}

\begin{lemma}
	The orthogonal projection of a positive semidefinite matrix onto a transpose-closed
	real matrix algebra is positive semidefinite.
\end{lemma}

\proof
This is a special case of Tomiyama's theorem, see \cite{Tom59}.\qed

Given an orthogonal basis for the Bose-Mesner algebra, we can compute orthogonal projections of matrices onto it---if $M$ is an $\binom{n}{k}\times\binom{n}{k}$ matrix, its orthogonal projection $\Psi(M)$ is
given by Gram-Schmidt:
\[
	\Psi(M) = \sum_i \frac{\ip{M}{A_i}}{\ip{A_i}{A_i}}A_i.
\]
Note that 
\[
	\ip{M-\Psi(M)}{A} = 0
\]
for any matrix $A$ in the Bose-Mesner algebra, and taking $A$ to be $J$ and $I$ in turn
yields that
\[
	\elsm{\Psi(M)} = \tr(M),\quad \tr(\Psi(M)) = \tr(M).
\]

We consider an example. For any family $\cF$ of $k$-subsets of $N$, we denote by $N_{\cF}$ the
matrix $xx^T$ where $x$ is the characteristic vector of $\cF$. Let $\cF$ be a $t$-intersecting family of $k$-subsets.
Then
\[
	\ip{A_r}{N_\cF} = \tr(AN_\cF) = x^TA_rx,
\]
which equals the number of pairs $(a\,\be)$ in $\cF\times\cF$ such that $|\al\cap\be|=k-r$.
Therefore $\ip{A_r}{N_\cF}=0$ if $r\ge k-t+1$.

\begin{lemma}
	Let $\cF$ be a $t$-intersecting family. Then $\Psi(N_\cF)$ is a positive semidefinite matrix
	lying in the span of $\seq A0{k-t}$ and
	\[
		\tr(\Psi(N_\cF)) = |\cF|, \quad \elsm(\Psi(N_\cF)) = |\cF|^2.
	\]
\end{lemma}

\proof
Observe that if $A$ lies in the Bose-Mesner algebra of the Johnson scheme, then
\[
	0 = \ip{M-\Psi(M)}{A} = \ip{M}{A} -\ip{\Psi(M)}{A},
\]
whence $\ip{\Psi(M)}{A}=\ip{M}{A}$. Therefore $\ip{\Psi(M)}{A_r}=\ip{M}{A_r}$, which proves
that $\Psi(N_\cF)$ lies in the span of $\seq A0{k-t}$. The remaining two claims follow
from the fact that $\Psi$ preserves trace.\qed

If we can show that the Bose-Mesner algebra of $J(n,k)$ contains a matrix $L$ such that:
\begin{enumerate}[(a)]
	\item 
	$L$ is positive semidefinite,
	\item
	$\Psi(N_\cF)\circ L = \ga I$ for some $\ga$,
	\item
	$\elsm(L)/\tr(L)=\binom{n}{t}/\binom{n-t}{k-t}$,
\end{enumerate}
then Lemma~\ref{lem:clqcoclq} implies that
\[
	|\cF| \le {\binom{n-t}{k-t}}.
\]
The key to Wilson's proof was to demonstrate that, provided
\[
	n\le(t+1)(k-t+1),
\] 
the matrix $I+\Om(n,k,t)$ satisfies these conditions.

Recall that a $t$-$(n,k,\la)$-design is a collection of subsets of size $k$ from an $n$-set such 
that any any subset of $t$ points from $V$ lies in exactly $\la$ blocks (aka $k$-sets).
If $\la=1$, we call the design a \textsl{Steiner system}. The construction of
Steiner systems for large $t$ is something of a mystery (to which we shall return),
but projective and affine planes of finite order provide examples with $t=2$ and
M\"obius planes give examples with $t=3$.  

\begin{lemma}
	Let $\cD$ be a $t$-$(n,k,1)$ design. Then $\Psi(N_\cD)$ is a positive semidefinite
	matrix lying in the span of $\seq A{k-t+1}{t}$ and
	\[
		\tr(\Psi(N_\cD)) = |\cD|,\quad \elsm(\Psi(N_\cD)) = |\cD|.
	\]
\end{lemma}

\begin{lemma}
	If a $t$-$(n,k,1)$-design exists, then a $t$-intersecting family
	of $k$-subsets of a set of size $v$ has size at most $\binom{n-t}{k-t}$.\qed
\end{lemma}

We can compute $\Psi(N_\cD)$ explicitly. If $\la_i$ denotes the number of blocks of $\cD$
that contain a given set of $i$ points and $0\le i\le t$, then
\[
	\la_i = \frac{\binom{n-i}{k-i}}{\binom{n-t}{k-t}}.
\]
If we define
\[
	\ga_s = \sum_{i=s}^t (-1)^{i-s}\binom{i}{s}\binom{k}{i}(\la_i-1),
\]
then from Exercise~1 in Chapter 8 of \cite{cg-km}, we find that
\[
	\Psi(N_\cD) = \sum_{s=0}^t \frac{\ga_s}{\binom{n-k}{k-s}\binom{k}{s}} A_{k-s}.
\]
 For $n,k,t$, we will denote by $M(n,k,t)$ the following:
\[
M(n,k,t) = \sum_{s=0}^t \frac{\ga_s}{\binom{n-k}{k-s}\binom{k}{s}} A_{k-s}. 
\]
If there exists a  $t$-$(n,k,1)$-design $\cD$ exists, then $\Psi(N_\cD) = M(n,k,t) $.
Observe that the matrix $M(n,k,t)$ is always well-defined, whether or not the design exists. 
We use Keevash's result \cite{Kee14} on the existence of $t$-designs to show that this
projection is equal to Wilson's matrix. For a second proof of Keevash's result, see \cite{GloKuhLoOs16}.

The following theorem is a restatement of Theorem 1.4 of \cite{Kee14} applied to $G=K_n^t$, 
in the language of block designs instead of hypergraphs. 

\begin{theorem}[Keevash]
	\label{thm:keevash} 
	For fixed $k$ and $t$, there exists $N$ such that for $n > N$, if $\binom{k-i}{t-i}$ 
	divides $\binom{n-i}{t-i}$ for $i = 0,\ldots,t-1$, then there exists a $t$-$(n,k,1)$ 
	block design.
\end{theorem}

We are now able to prove the following.

\begin{theorem} 
	For any $n\geq k \geq t$, we have that  $M(n,k,t) = \Omega(n,k,t) + I$. 
\end{theorem}

\proof 
Fix $k$ and  $t$. Let \[f_r(n) = \theta_r(M(n,k,t))\] and \[g_r(n) = \theta_r(\Omega(n,k,t)) + 1.\] 
If there exists a  $t$-$(n,k,1)$-design $\cD$ exists, then $f_r(n) = g_r(n)$ 
for $r = 0,\ldots, t$. By Theorem~\ref{thm:keevash}, we have that $f_r(n)$ and $g_r(n)$ are equal 
for infinitely many $n$. Consider $h_r(n) = f_r(n) - g_r(n)$. We see that $h_r(n)$ is a rational 
function whose numerator $p(n)$ is a polynomial. Since $p(n) =0 $ infinitely often, we have that $p(n) = 0$
 and so $h_r(n) =0$. We thus have that  $f_r(n)= g_r(n)$ for all $n$. This shows that  
$M(n,k,t) = \Omega(n,k,t) + I$ for all $n$.\qed 

%

\begin{thebibliography}{1}

\bibitem{AhlKha96}
R.~Ahlswede and L.H. Khachatrian.
\newblock The complete nontrivial-intersection theorem for systems of finite
  sets.
\newblock {\em Journal of Combinatorial Theory. Series A}, 76(1):121--138,
  1996.

\bibitem{AhlKha97}
R.~Ahlswede and L.H. Khachatrian.
\newblock The complete intersection theorem for systems of finite sets.
\newblock {\em European Journal of Combinatorics}, 18(2):125--136, 1997.

\bibitem{GloKuhLoOs16}
S.~{Glock}, D.~{K{\"u}hn}, A.~{Lo}, and D.~{Osthus}.
\newblock {The existence of designs via iterative absorption}.
\newblock {\em ArXiv e-prints}, November 2016.

\bibitem{cg-km}
Chris Godsil and Karen Meagher.
\newblock {\em {Erd\H{o}s}-{K}o-{R}ado Theorems: Algebraic Approaches}, volume
  149 of {\em Cambridge Studies in Advanced Mathematics}.
\newblock Cambridge University Press, Cambridge, 2016.

\bibitem{Kee14}
P.~{Keevash}.
\newblock {The existence of designs}.
\newblock {\em ArXiv e-prints}, January 2014.

\bibitem{Tom59}
Jun Tomiyama.
\newblock On the projection of norm one in {$W^*$-algebras}, {III}.
\newblock {\em Tohoku Math. J. (2)}, 11(1):125--129, 1959.

\bibitem{Wilson1984}
R.M. Wilson.
\newblock The exact bound in the {Erd\H{o}s}-{K}o-{R}ado theorem.
\newblock {\em Combinatorica}, 4(2-3):247--257, 1984.

\end{thebibliography}

\end{document}